\documentclass[12pt]{article}

\newtheorem{Theorem}{Theorem}[section]
\newtheorem{Proposition}[Theorem]{Proposition}
\newtheorem{Lemma}[Theorem]{Lemma}

\newtheorem{Corollary}[Theorem]{Corollary}

\def\buildreb#1\over#2{\mathrel{\mathop{\kern0pt #1}\limits_{#2}}}

\def\fin{\end{document}}

\newtoks\baspagetitre
\baspagetitre={\hfill{}}

\newif\ifpagetitre \pagetitretrue 
\def\footnoterule{\kern -3pt\hrule width 2truein \kern 2.4pt}
\renewcommand{\thepage}{\ifpagetitre\the\baspagetitre\global\pagetitrefalse
 \else \arabic{page}\fi}

\begin{document}
 
 \title{\sc Simultaneous Diophantine approximation with a divisibility condition}
 
  \author{\bf Bernard de Mathan}
 
  \date{}
\maketitle
{\vskip3cm}

\leftskip = 0pt
\rightskip = 2pt
{\bf R\'esum\'e}. Dans un article pr\'ec\'edent (\cite{Mon}), nous \'etudiions des approximations rationnelles simultan\'ees dans ${\bf R}^2$ qui pr\'esentent une certaine analogie avec les fractions continues. Nous obtenions des r\'esultats autour de la conjecture de Littlewood en utilisant de telles approximations. Nous montrons ici que ces r\'esultats restent vrais si l'on ajoute des conditions de divisibilit\'e.\par\ \par
{\bf Abstract.} In a previous paper (\cite{Mon}), we studied certain sequences of simultaneous rational approximations in ${\bf R}^2$ which present some analogy with the continued fractions. We got results around the Littlewood conjecture by using such approximations. Here we show that these results also hold when we add divisibility conditions.\par\ \par
{\bf Keywords} Littlewood conjecture . Simultaneous Diophantine approximation . Divisibility\par\ \par
\section{Introduction.} 
 The Littlewood  conjecture in simultaneous diophantine approximation claims that for every pair $(\alpha,\beta)$ of real numbers, there exists for every $\epsilon>0$, a triple of integers $(q,r,s)$ with $q>0$, such that
$$q\vert q\alpha-r\vert\vert q\beta-s\vert\le\epsilon.\eqno\hbox{(L)}$$
Some results are known about this conjecture. Cassels and Swinnerton-Dyer proved that every pair of numbers in a cubic field satisfies this statement \cite{CS}.
A more precise result was given by Peck \cite{Pe}.  Einsiedler, Katok and Lindenstrauss proved that the set of exceptional pairs $(\alpha,\beta)$ is very small, that is to say that it has Hausdorff dimension zero \cite{EKL}.

Here we consider other conditions. For instance, we may ask whether for any pair $(\alpha,\beta)\in{\bf R}^2$, there is a positive constant $C=C(\alpha,\beta)$ such that for any positive number $D$ there exist integers $q$, $r$, and $s$, with $q>0$ and  $$q\vert q\alpha-r\vert\vert q\beta-s\vert\le C, \qquad D|q\eqno\hbox{(LB)}$$ (the notation $D|q$ means that $D$ divides $q$). This condition is weaker than the Littlewood conjecture since if a triple $(q,r,s)$ satisfies (L) with
$\epsilon=D^{-3}$, then the triple $(Dq,Dr,Ds)$ satisfies (LB) with $C=1$. We do not know whether (LB) is satisfied for any pair $(\alpha,\beta)\in{\bf R}^2$.

Here we consider a stronger condition. For any pair $(\alpha,\beta)\in{\bf R}^2$, we ask whether there exists for every $\epsilon>0$ and every integer $D\ge1$, a triple of integers $(q,r,s)$ with $q>0$ such that
$$q|q\alpha-r||q\beta-s|\le\epsilon,\qquad D|q,\qquad{\rm g.c.d.}(q,r,s)=1\eqno\hbox{(Ldiv)}$$ 
We can again consider the following stronger condition:  for any $\epsilon>0$ and for any positive integer $D$, does exist a triple of integers $(q,r,s)$, with $q>0$, such that
$$q|q\alpha-r||q\beta-s|\le\epsilon,\qquad q\equiv r\equiv0\ {\rm mod}\ D,\qquad {\rm g.c.d.}(q,r,s)=1\ ?\eqno\hbox{(Ldiv2)}$$ Frequently we simply denote the greatest common divisor ${\rm g.c.d.}(q,r,s)$ by $(q,r,s)$.\par
Notice that in order to get (Ldiv2), it is enough to require that $q\equiv r\equiv0\ {\rm mod}\ D$, and $(D,s)=1$. Indeed, if we consider $\Delta={\rm g.c.d.}(q,r,s)$, we then have $(D,\Delta)=1$ and if we put $q_1=q/\Delta$, $r_1=r/\Delta$ and $s_1=s/\Delta$, then we have
$$q_1|q_1\alpha-r_1||q_1\beta-s_1|\le{\epsilon\over\Delta^3}\le\epsilon,\qquad q_1\equiv r_1\equiv0\ {\rm mod}\ D,\qquad {\rm g.c.d.}(q_1,r_1,s_1)=1,$$ which is (Ldiv2).
\par If $\alpha$ is a rational number or has an infinite Markov constant, then condition (Ldiv2) is satisfied for any $\beta\in{\bf R}$. Indeed, in this case, for every $\epsilon>0$ and every positive integer $D$, there exist integers $q$, $r$ and $s$, $q>0$, such that $$q|q\alpha-r|\le{\epsilon\over D^2(D+1)},\qquad |q\beta-s|<1.
$$ Taking $q'=Dq$, $r'=Dr$, and choosing an integer $s'$ such that $(D,s')=1$ and $|q'\beta-s'|<D+1$ (for instance $s'=Ds+1$), we get
$$q'|q'\alpha-r'||q'\beta-s'|\le\epsilon,\qquad q'\equiv r'\equiv0\ {\rm mod}\ D, \qquad{\rm g.c.d.}(D,s')=1,$$ which implies (Ldiv2) by the above remark.\par
The converse is true when $\beta$ is a rational number, 
that is to say that, given $\beta\in{\bf Q}$, the pair $(\alpha,\beta)$ satisfies (Ldiv2) if {\sl and only if} $\alpha$ has an infinite Markov constant (or is a rational number). Indeed, suppose that $\beta=B/C$, where  $B$ and $C$ are integers with $(B,C)=1$ and $C>0$. In order to get (Ldiv2) with a triple $(q,r,s)$, we must take $s$ with $(s,D)=1$. Hence if $q$ is a multiple of $D$, $q=DQ$, then
 $$|DQ\beta-s|={|DBQ-sC|\over C},$$ thus if $D\not|\ C$, we have then $DBQ-sC\ne0$ and $$|DQ\beta-s|\ge{1\over C}\cdot$$ Accordingly if we have (Ldiv2) with $r=DR$, we must have
 $$Q|Q\alpha-R|\le{C\over D^2}\epsilon,$$ thus $\alpha$ has an infinite Markov constant. Note that the condition (Ldiv2) is not symmetrical since for instance, if $\beta$ is a rational number, then the pair $(\beta,\alpha)$ satisfies (Ldiv2) for every $\alpha\in{bf R}$ while $(\alpha,\beta)$ satisfies (Ldiv2) only if $\alpha$ has an infinite Markov constant or is a rational number.\par
 We do not know any example of a pair $(\alpha,\beta)\in{\bf R}^2$ which does not satisfy (Ldiv). For (Ldiv2), we do not know any exception other than the pairs $(\alpha,\beta)$ where $\alpha$ is a real number with a finite Markov constant and $\beta$ is a rational number.\par
 We shall obtain some results about these conditions by using the notion of recurrent word. This notion was already used in view of another version of the Littlewood conjecture in \cite{Bu}.
The results which we shall give involve pairs $(\alpha,\beta)$ of real numbers which are simultaneously  badly approximable, namely for which there exists a real constant $C=C(\alpha,\beta)>0$ such that 
$$q^{1/2}\max\{\vert q\alpha-r\vert, \vert q\beta-s\vert\}\ge C\eqno\hbox{(Bad2)}$$
for every triple of integers $(q,r,s)$ with $q>0$.\par We use the Vinogradov notations $\ll$, $\gg$ and $\asymp$ (we write $A\ll B$ to mean that $0\le A\le KB$ for some positive constant $K$; $A\asymp B$ means that $A\ll B$ and $B\ll A$). The dependence of the implied constants will be generally clear.\par We start from simultaneous approximations which we have studied in \cite{Mon}, where we proved the following result:
\begin{Proposition} Let $(\alpha,\beta)$ be a pair of real numbers. If $(\alpha,\beta)$ satisfies {\rm (Bad2)}, then there exists a sequence of triples of integers $(q_n,r_n,s_n)_{n\ge0}$, with $q_n>0$, satisfying the following conditions:
$$
q_n^{1/2}|q_n\alpha-r_n|\ll1,\ \ \ \ \ q_n^{1/2}|q_n\beta-s_n|\ll1,   \eqno\hbox{\rm(1.1)}
$$
$$
q_n\asymp q_{n-1},      \hspace{3cm}     (n\ge1);       \eqno\hbox{\rm(1.2)}
$$ there exists a positive integer $\chi$ and a constant $K>1$ such that
$$q_{n+\chi}\ge Kq_n{\hskip2cm} (n\ge0),                 \eqno\hbox{\rm(1.3)}
$$ 
$$
\left\vert\matrix{q_n&r_n&s_n\cr q_{n-1}&r_{n-1}&s_{n-1} \cr  q_{n-2}&r_{n-2}&s_{n-2}}\right\vert\ne0 {\hskip2cm} (n\ge2).                       \eqno\hbox{\rm(1.4)}$$ 
Conversely, if there exists a sequence of triples of integers $(q_n,r_n,s_n)$ with $q_n$ positive and unbounded, such that conditions {\rm(1.1)}, {\rm(1.2)} and {\rm(1.4)} are satisfied, then the pair $(\alpha,\beta)$ satisfies {\rm (Bad2)}.
\end{Proposition} 
We also proved in \cite{Mon}:
\begin{Lemma} Let $(\alpha,\beta)$ be a pair of real numbers. Suppose that $(q_n,r_n,s_n)_{n\ge0}$ is a sequence of triples of integers, with $q_n$ positive, satisfying the conditions {\rm(1.1)}, {\rm(1.2)}, {\rm(1.4)}. Then there exist sequences $(a_n)_{n\ge3}$, $(b_n)_{n\ge3}$ and $(c_n)_{n\ge3}$  of rational numbers of bounded height such that for any $n\ge2$, we have
$$q_{n+1}=a_{n+1}q_n+b_{n+1}q_{n-1}+c_{n+1}q_{n-2},    \eqno\hbox{\rm(1.5)}$$
$$r_{n+1}=a_{n+1}r_n+b_{n+1}r_{n-1}+c_{n+1}r_{n-2},      \eqno\hbox{\rm(1.6)}$$ and
$$s_{n+1}=a_{n+1}s_n+b_{n+1}s_{n-1}+c_{n+1}s_{n-2}.          \eqno\hbox{\rm(1.7)}$$
\end{Lemma}
We shall use these properties to state our results. Recall again the definition:\par
{\bf Definition}. {\sl An infinite word $(\mu_n)_{n\ge0}$, where the $\mu_n$ lie in a finite set - the alphabet -, is recurrent if for each integer $N\ge0$, there exists an integer $k>0$ such that:}
$$\mu_n=\mu_{n+k}\hspace{3cm}(0\le n\le N).$$
We can then state our results:
\begin{Theorem} Let $(\alpha,\beta)$ be a pair of real numbers satisfying {\rm (Bad2)}. Let $(q_n,r_n,s_n)$ be a sequence of triples of integers, with $q_n>0$, satisfying conditions {\rm(1.1)}, {\rm(1.2)}, {\rm(1.3)} and {\rm(1.4)}. Suppose that the sequence of the triples
of coefficients $(a_n,b_n,c_n)_{n\ge3}$ satisfying {\rm(1.5)}, {\rm (1.6)} and {\rm (1.7)} is an infinite recurrent word. Then the pair $(\alpha,\beta)$ satisfies the condition {\rm (Ldiv2)}. More precisely, for any $\epsilon>0$ and any integer $D>0$, there exist integers $Q$, $R$, and $S$, with $Q>0$ and ${\rm g.c.d.}(Q,R,S)=1$, such that
$$ 
\vert Q\alpha-R\vert\le\epsilon Q^{-1/2},\qquad\vert Q\beta-S\vert\ll_{\alpha,\beta} Q^{-1/2},\qquad Q\equiv R\equiv0\ {\rm mod}\ D\eqno\hbox{\rm(1.8)}$$ (the constant which is implied in {\rm(1.8)} only depends upon $\alpha$ and $\beta$).
\end{Theorem} Note that the condition that the word $(a_n,b_n,c_n)_{n\ge3}$ is recurrent makes sense since this sequence takes only a finite number of values (which are triples of rational numbers of bounded height). 
\par
In \cite{Mon}, we had already proved that under the hypotheses of Theorem 1.3, the pair $(\alpha,\beta)$ satisfies the Littlewood conjecture (L).\par
Theorem 1.3 asserts in particular that there is a constant $C_0=C_0(\alpha,\beta)>0$ such that for each positive integer $D$, there exist infinitely many positive integers $Q$ satisfying
$$Q^{1/2}\max\{\Vert Q\alpha\Vert,\Vert Q\beta\Vert\}\le C_0,\qquad \ \qquad
Q\equiv0\ {\rm mod}\ D,\eqno\hbox{(1.9)}$$ (where $\Vert.\Vert$ denotes the distance to the nearest integer).\par In \cite{Bg}, Bugeaud studies the Lagrange constant of a pair $(\alpha,\beta)$ of real numbers:
$$c(\alpha,\beta)={\buildreb{\lim\inf}\over{Q>0}}Q^{1/2}\max\{\Vert Q\alpha\Vert,\Vert Q\beta\Vert\}.$$ Here inequality (1.9) implies that:
$$c(D\alpha,D\beta)\le C_0D^{-1/2}.\eqno\hbox{(1.10)}$$ Hence we have got:
\begin{Corollary} Under the hypotheses of Theorem {\rm1.3}, we have for each positive integer $D$
$$c(D\alpha,D\beta)\ll_{\alpha,\beta} D^{-1/2},\eqno\hbox{\rm(1.10)}$$ where the implied constant only depends upon the pair $(\alpha,\beta)$. In particular
$$\lim_{D\to\infty}c(D\alpha,D\beta)=0.$$ \end{Corollary}
Corollary 1.4 applies in particular when $(1,\alpha,\beta)$ is a basis of a cubic field since in \cite{Mon} we have proved that, in this case, there is a sequence of triples of integers $(q_n,r_n,s_n)$ with $q_n>0$, satisfying (1.1), (1.2), (1.3) and (1.4), for which the sequence $(a_n,b_n,c_n)$ defined by conditions (1.5), (1.6) and (1.7), is recurrent. However in the case of
a pair $(\alpha,\beta)$ of numbers in the same cubic field, inequality (1.10) was already proved by Bugeaud in \cite{Bg}. We do not know whether (1.10) is true for any pair of simultaneously  badly approximable numbers. In \cite{Bg}, Bugeaud states some problems about the weaker condition $$\sup_{D\ge1} c(D\alpha,D\beta)<+\infty.\eqno\hbox{(1.10')}$$ 

In the case where $(1,\alpha,\beta)$ is a basis of a cubic field $E$, Theorem 1.3 also applies, however (Ldiv) was already ensured by an argument of \cite{dMT} (Th\'eor\`eme 3.1). Here we shall give a more effective version by using the classical method of Peck \cite{Pe}. Let us define a multiplicative function $\psi$ on the set 
${\bf Z}_{>0}$ of the positive integers such that for each prime number $p$ and each positive integer $\nu$, $$\psi(p^\nu)=(p^6-1)p^{\nu}.$$ That means that for every positive integer $D$, we have $$\psi(D)=D\prod_{p|D}(p^6-1).$$ 
We shall prove:
\begin{Theorem} Let $E$ be a cubic field, and let $\alpha$ and $\beta$ be numbers such that $(1,\alpha,\beta)$ is a linear basis of $E$.\par
{\rm a)} If $E$ has a unique real embedding, then for every $0<\epsilon\le1$ and every positive integer $D$, there exists a triple $(Q,R,S)$ of integers with $Q>1$ and ${\rm g.c.d.}(Q,R,S)=1$, satisfying $$ 
\vert Q\alpha-R\vert\le\epsilon Q^{-1/2},\qquad\vert Q\beta-S\vert\ll_{\alpha,\beta} Q^{-1/2},\qquad Q\equiv R\equiv0\ {\rm mod}\ D,\eqno\hbox{\rm(1.8)}$$ with
$$\log Q\ll_{\alpha,\beta}\epsilon^{-1}\psi(D).\eqno\hbox{\rm(1.13)}$$
{\rm b)} If the field $E$ is totally real, then the same holds when replacing {\rm(1.13)} by
$$\log Q\ll_{\alpha,\beta}\epsilon^{-1}\psi^2(D).\eqno\hbox{\rm(1.14)}$$  
(All the implied constants only depend upon $\alpha$ and $\beta$).\end{Theorem}
In \cite{Bg}, Theorem 1.3, Bugeaud claims that under the hypotheses of Theorem 1.5, for any prime number $p$,
$${\buildreb{\lim\inf}\over{Q>0}}Q^{1/2}(\log Q)|Q|_p\max\{\Vert Q\alpha\Vert,\Vert Q\beta\Vert\}<+\infty.\eqno\hbox{(1.7)[4]}$$
This is implied by our Theorem 1.5 in the case a). Indeed, if we take $\epsilon=1$ and $D=p^\nu$, where $\nu$ is a positive integer, then we deduce from Theorem 1.5 that there is a positive integer $Q$ such that
 $$ 
\vert Q\alpha-R\vert\le Q^{-1/2},\qquad\vert Q\beta-S\vert\ll_{\alpha,\beta} Q^{-1/2},\qquad Q\equiv 0\ {\rm mod}\ p^\nu$$ with $$0<\log Q\ll_{\alpha,\beta}(p^6-1)p^\nu,$$ i.e.,
$$|Q|_p\log Q\ll_{\alpha,\beta,p}1.$$ Then (1.8) leads to (1.7)[4].
 However it seems that the proof of (1.7)\cite{Bg} which is given in \cite{Bg} works only in the case a), because for any unit $\eta$ of $E$, and any complex not real embedding $\sigma$ of $E$, we have then $|\sigma(\eta)|=|\eta|^{-1/2}$. This condition fails in the case b), it is only possible to construct ``dominant'' units $\eta$ such that for any real embedding $\sigma$ other than identity, $|\sigma(\eta)|$ is ``{\sl nearly}'' equal to $|\eta|^{-1/2}$. For this reason, it seems that the proof of (1.7)\cite{Bg} is incomplete in the case b). In this case, Theorem 1.5 leads to the inequality: $${\buildreb{\lim\inf}\over{Q>0}}Q^{1/2}(\log Q)^{1/2}|Q|_p\max\{\Vert Q\alpha\Vert,\Vert Q\beta\}<+\infty,$$ which can be already found in \cite{dMT}.

 \section{Proof of Theorem 1.3}
For convenience, we reproduce here some preliminaries which are analogous to those which we used in \cite{Mon}.
Let $q_n$, $r_n$, $s_n$, be sequences of integers with $q_n>0$, and suppose that conditions (1.1), (1.2), (1.3) and (1.4) are satisfied. For $2\le m<k$, consider the determinants:
$$ x_{m,k}=\left\vert\matrix{q_m&r_{m}&q_k\cr  q_{m-1}&r_{m-1} & q_{k-1} \cr q_{m-2}&r_{m-2}&q_{k-2}}\right\vert, \ \ \ y_{m,k}=\left\vert\matrix{q_m&r_{m}&r_k\cr  q_{m-1}&r_{m-1} & r_{k-1} \cr q_{m-2}&r_{m-2}&r_{k-2}}\right\vert,
\ \ \ z_{m,k}=\left\vert\matrix{q_m&r_{m}&s_k\cr  q_{m-1}&r_{m-1} & s_{k-1} \cr q_{m-2}&r_{m-2}&s_{k-2}}\right\vert.\eqno\hbox{(2.1)}
$$ These determinants are integers.
\begin{Lemma}
We have the estimations:
$$
\max\{|x_{m,k}|, |y_{m,k}|, |z_{m,k}|\}\ll q_m^{1/2}q_k, \eqno\hbox{\rm(2.2)}$$
 and      $$
\max\{|x_{m,k}\alpha -y_{m,k}|, |x_{m,k}\beta-z_{m,k}|\}\ll q_m^{1/2}q_k^{-1/2}. \eqno\hbox{\rm(2.3)}
$$
Furthermore there exists a constant $L>0$, which depends only upon $\alpha$ and $\beta$, such that if $k-m\ge L$, then we have $x_{m,k}\ne0$.
\end{Lemma}   {\sl Proof}. Set
$$q_n\alpha-r_n=\alpha_n, {\hskip3cm} q_n\beta-s_n=\beta_n.$$
We get inequality (2.2) by writing
$$x_{m,k}=-\left\vert\matrix{q_m&\alpha_{m}&q_k\cr  q_{m-1}&\alpha_{m-1} & q_{k-1} \cr q_{m-2}&\alpha_{m-2}&q_{k-2}}\right\vert,\ \ y_{m,k}=-\left\vert\matrix{q_m&\alpha_{m}&r_k\cr  q_{m-1}&\alpha_{m-1} & r_{k-1} \cr q_{m-2}&\alpha_{m-2}&r_{k-2}}\right\vert,\ \
z_{m,k}=-\left\vert\matrix{q_m&\alpha_{m}&s_k\cr  q_{m-1}&\alpha_{m-1} & s_{k-1} \cr q_{m-2}&\alpha_{m-2}&s_{k-2}}\right\vert.
   \eqno\hbox{(2.4)}$$ Then we get (2.2) by using (1.1) and (1.2).
 It follows from (2.4) that 
$$ x_{m,k}\alpha-y_{m,k}=-\left\vert\matrix{q_m&\alpha_{m}&\alpha_k\cr  q_{m-1}&\alpha_{m-1} & \alpha_{k-1} \cr q_{m-2}&\alpha_{m-2}&\alpha_{k-2}}\right\vert,\qquad x_{m,k}\beta-z_{m,k}=-\left\vert\matrix{q_m&\alpha_{m}&\beta_k\cr  q_{m-1}&\alpha_{m-1} & \beta_{k-1} \cr q_{m-2}&\alpha_{m-2}&\beta_{k-2}}\right\vert.
\eqno\hbox{(2.5)}$$ Then we get
$$\max\{|x_{m,k}\alpha -y_{m,k}|, |x_{m,k}\beta-z_{m,k}|\}\ll q_m^{1/2}q_k^{-1/2},$$ which is (2.3).\par
Let us prove that $x_{m,k}\ne0$. First note that it is not possible that
$x_{m,k}=y_{m,k}=z_{m,k}=0$. Indeed set $$\left\vert\matrix{q_n&r_{n}&s_n\cr  q_{n-1}&r_{n-1} & s_{n-1} \cr q_{n-2}&r_{n-2}&s_{n-2}}\right\vert=\Delta_n.$$ By (1.4), we have $\Delta_n\ne0$ for each $n$. As $\Delta_m\ne0$, if $x_{m,k}=y_{m,k}=z_{m,k}=0$, then the vectors $(q_k,q_{k-1},q_{k-2})$, $(r_k,r_{k-1},r_{k-2})$, and $(s_k,s_{k-1},s_{k-2})$ should belong to the plane ${\bf Q}(q_m,q_{m-1},q_{m-2})+{\bf Q}(r_m,r_{m-1},r_{m-2})$, contrarily to the fact that $\Delta_k\ne0$. Now notice that by (1.2) and (1.3), we have $$q_k\gg K^{k-m\over \chi}q_m, \eqno\hbox{(1.3')}$$ 
thus we deduce from (2.3) that $$\max\{|x_{m,k}\alpha -y_{m,k}|, |x_{m,k}\beta-z_{m,k}|\}\ll K^{m-k\over2\chi}.$$ 
Consequently there exists an integer $L>0$, only depending upon $\alpha$ and $\beta$, such that, for $k-m\ge L$, we have 
$$\max\{|x_{m,k}\alpha -y_{m,k}|, |x_{m,k}\beta-z_{m,k}|\}<1.$$  Hence, if $x_{m,k}=0$, then we must have $y_{m,k}=z_{m,k}=0$, which is impossible. Thus we get $x_{m,k}\ne0$.\par
It will be useful in the sequel to note that $$|\Delta_n|\asymp1.\eqno\hbox{(2.6)}$$ Indeed, $\Delta_n$ is a non zero integer, and as we can write $$
\Delta_n=\left\vert\matrix{q_n&\alpha_{n}&\beta_n\cr  q_{m-1}&\alpha_{m-1} & \beta_{m-1} \cr q_{m-2}&\alpha_{m-2}&\beta_{m-2}}\right\vert,$$ we have
$$1\le|\Delta_n|\ll1.$$  

We regard now the coefficients $a_n$, $b_n$, and $c_n$, defined by (1.5), (1.6) and (1.7). Put $$\mu_n=(a_n,b_n,c_n).$$ The crucial point is the following:
\begin{Lemma} Let $m$ and $k$ be integers with $2<m<k$.
              Suppose that for each integer $i$ with $0\le i\le m-3$, we have 
$$         \mu_{k-i}=\mu_{m-i}   ,     $$ i.e.,
$$ \mu_k=\mu_m,\hspace {2cm}\dots, \hspace{2cm} \mu_{k-m+3}=\mu_{3}.$$
Then            $$
{x_{m,k}\over\Delta_m}={x_{2,k-m+2}\over\Delta_2}, {\hskip1.5cm}
{y_{m,k}\over\Delta_m}={y_{2,k-m+2}\over\Delta_2},{\hskip1.5cm}
{z_{m,k}\over\Delta_m}={z_{2,k-m+2}\over\Delta_2}\cdot
\eqno\hbox{\rm(2.7)}$$
Moreover, we have $$q_k\asymp q_mq_{k-m}\eqno\hbox{\rm(2.8)}$$  (the constants which are implied in {\rm(2.8)} only depend upon the pair $(\alpha,\beta)$ and the sequence of triples $(q_n,r_n,s_n)$ chosen).
\end{Lemma}
{Proof}. We introduce the matrices
 $$A_n=\left(\matrix{a_n&b_n&c_n\cr1&0&0\cr0&1&0}\right)$$ and 
 $$Q_n=\left(\matrix{q_n&r_{n}&s_n\cr  q_{n-1}&r_{n-1} & s_{n-1} \cr q_{n-2}&r_{n-2}&s_{n-2}}\right).$$  
Equations (1.5), (1.6) and (1.7) can be written
$$
Q_n=A_nQ_{n-1} \ \ \ \ \ (n\ge3).
$$ Since $A_m=A_k$, we have
$$\left(\matrix{q_m&r_{m}&q_k\cr  q_{m-1}&r_{m-1} & q_{k-1} \cr q_{m-2}&r_{m-2}&q_{k-2}}\right)= A_m\left(\matrix{q_{m-1}&r_{m-1}&q_{k-1}\cr  q_{m-2}&r_{m-2} & q_{k-2} \cr q_{m-3}&r_{m-3}&q_{k-3}}\right).$$   
Therefore we get $$x_{m,k}=(\det A_m)x_{m-1,k-1}.$$ Now, as
$$\Delta_m=(\det A_m)\Delta_{m-1},$$ we get
$${x_{m,k}\over\Delta_m}={x_{m-1,k-1}\over\Delta_{m-1}}\cdot$$
By repeating this argument, we get$$
{x_{m,k}\over\Delta_m}={x_{2,k-m+2}\over\Delta_2}\cdot$$ The same proof applies to the equalities $${y_{m,k}\over\Delta_m}={y_{2,k-m+2}\over\Delta_2},{\hskip1.5cm}
{z_{m,k}\over\Delta_m}={z_{2,k-m+2}\over\Delta_2},$$ and thus (2.7) is proved.\par
We prove then (2.8). Setting $k-m=h$, by (1.4) and (1.2), we can find rational numbers $U$, $V$, $W$, such that we have for $j=0,1,2,$
$$q_{j+h}=Uq_j+Vr_j+Ws_j,\eqno\hbox{(2.9)}$$ and
$$\max\{|U|,|V|,|W|\}\ll q_h.$$ Now as $(a_{n+h},b_{n+h},c_{n+h})=(a_n,b_n,c_n)$ for each integer $n$ with $3\le n\le m$, equality (2.9) holds true for every $j$ with $0\le j\le m$.
Putting $$U'=U+V\alpha+W\beta,$$  we can write
$$q_{j+h}=U'q_j-V\alpha_j-W\beta_j\hspace{3cm}(0\le j\le m).\eqno\hbox{(2.10)}$$
We have $|V\alpha_j+W\beta_j|\ll q_j^{-1/2}q_h$. Let us write this inequality in an explicit form: $$|V\alpha_j+W\beta_j|\le K_1q_j^{-1/2}q_h$$
with a constant $K_1>0$, depending only upon the sequence $(q_n,r_n,s_n)$. We need also to write (1.3') in an explicit form:  there is a constant $K_2>0$ such that 
for any pair of integers $(n,l)$ with $n\ge l\ge0$, we have $$q_n\ge K_2K^{n-l\over\chi}q_l.\eqno\hbox{(1.3'')}.$$ Hence $q_j\ge K_2K^{j/\chi}q_0$ and
$q_{h+j}\ge K_2K^{j/\chi}q_h$, thus we get 
 $$|V\alpha_j+W\beta_j|\le K_1K_2^{-{1\over2}}K^{-{j\over2\chi}}q_0^{-\frac12}q_h\le K_1K_2^{-{3\over2}}K^{-{3j\over2\chi}}q_0^{-\frac12}q_{h+j}.$$ Since $K>1$, we can find a positive integer $j_0$, which only depends upon $(\alpha,\beta)$, such that $K^{3j_0/(2\chi)}\ge2K_1K_2^{-3/2}q_0^{-1/2}$. If $2<m<j_0$, as $j_0$ is a constant, we then have by (1.2), $q_m\asymp1$ and $q_{k-m}\asymp q_k$, thus (2.8) is trivial. If $m\ge j_0$, then we deduce from the above inequality that for any $j$ such that $j_0\le j\le m$, we have
 $$|V\alpha_{j}+W\beta_{j}|\le\frac{q_{h+j}}2\cdot$$ If we compare this inequality with (2.10), then we get 
 $$ \frac12q_{h+j}\le U'q_{j}\le \frac32q_{h+j}\hspace{3cm}(j_0\le j\le m).$$
We take successively $j=m$ and $j=j_0$ in this inequality: for $j=m$, we get
 $$U'q_m\asymp q_k,$$ and for $j=j_0$ (a constant, so that by (1.2), $q_{h+j_0}\asymp q_h$ and $q_{j_0}\asymp1$):
 $$U'\asymp q_h,$$ Thus we have
 $$q_k\asymp q_hq_m,$$ which is (2.8).\par
 We can now achieve  the proof of Theorem 1.3.
 Since the word $(\mu_n)_{n\ge3}$ is recurrent, there exists a strictly increasing sequence of integers $(k_i)$ such that
 $$\mu_{k_{i+1}}=\mu_{k_i}, \ \cdots , \ \  \mu_{k_{i+1}-k_i+3}=\mu_{3}.\eqno\hbox{(2.11)}$$
Indeed, we may construct inductively such a sequence by starting for instance from $k_0=3$, and given $k_i$, we choose $k_{i+1}>k_i$ such that the finite word $(\mu_k)_{3\le k\le k_i}$  is identical to the word $(\mu_k)_{k_{i+1}-k_i+3\le k\le k_{i+1}}$. Note that for any pair of integers
 $(i,h)$ with $h> i\ge0$, we have 
$$
\mu_{k_h-j}=\mu_{k_i-j} \hspace{3cm} (0\le j\le k_i-3).\eqno\hbox{(2.11')}$$
Replacing the sequence $(k_i)$ by a subsequence, we can also suppose that for each pair of integers $(i,h)$ with $h> i\ge0$, we have $k_h\ge k_i+L$, where the constant $L$ is defined as in Lemma 2.1. 
By Lemma 2.2, we have $$
x_{k_i,k_h}={\Delta_{k_i}\over\Delta_2}x_{2,k_h-k_i+2},$$
and by (2.6) and Lemma 2.1 with $m=2$ and $k=k_h-k_i+2$,
$$0<|x_{k_i,k_h}|\ll q_{k_h-k_i} .\eqno\hbox{\rm(2.12)}$$
By Lemma 2.2 and inequality (2.6), we have $|x_{k_i,k_h}\beta-z_{k_i,k_h}|\asymp
|x_{2,k_h-k_i+2}\beta-z_{2,k_h-k_i+2}|$, thus, by Lemma 2.1,
$$|x_{k_i,k_h}\beta-z_{k_i,k_h}|\ll q_{k_h-k_i}^{-1/2}.\eqno\hbox{(2.13)}$$
Now there exists an infinite set $H$ of non negative integers such that for every pair $(i,h)\in H\times H$ with $0\le i<h$, we have $$q_{k_i-j}\equiv q_{k_h-j}\ \ {\rm mod}\ D,\qquad r_{k_i-j}\equiv r_{k_h-j}\ \ {\rm mod}\ D\qquad(j=0,1,2),$$ which implies by (2.1) that $$x_{k_ i,k_h}\equiv y_{k_ i,k_h}\equiv0\ \ {\rm mod}\ D.\eqno\hbox{(2.14)}$$ Using the compactness of the projective space ${\bf P}_2({\bf R})$, we see that there are pairs of integers $(i,h)\in H\times H$ with $h> i\ge0$ for which there is a real number $\lambda_{i,h}$ such that
$$|\alpha_{k_{h}-j}-\lambda_{i,h}\alpha_{k_i-j}|\le\epsilon\max_{0\le j\le 2}|\alpha_{k_{h}-j}|\ll\epsilon q_{k_h}^{-1/2} \hspace{2cm}(j=0,1,2).$$ As by (2.5)
$$ x_{k_i,k_h}\alpha-y_{k_i,k_h}=-\left\vert\matrix{q_{k_i}&\alpha_{k_i}&\alpha_{k_h}\cr  q_{k_i-1}&\alpha_{k_i-1} & \alpha_{k_h-1} \cr q_{k_i-2}&\alpha_{k_i-2}&\alpha_{k_h-2}}\right\vert =-\left\vert\matrix{q_{k_i}&\alpha_{k_i}&\alpha_{k_h}-\lambda_{i,h}\alpha_{k_i}\cr  q_{k_i-1}&\alpha_{k_i-1} & \alpha_{k_h-1}-\lambda_{i,h}\alpha_{k_i-1} \cr q_{k_i-2}&\alpha_{k_i-2}&\alpha_{k_h-2}-\lambda_{i,h}\alpha_{k_i-2}}\right\vert,
 $$
we get
$$|x_{k_i,k_h}\alpha-y_{k_i,k_h}|\ll\epsilon q_{k_i}^{1/2}q_{k_h}^{-1/2},$$ hence by (2.8) :
$$|x_{k_i,k_h}\alpha-y_{k_i,k_h}|\ll \epsilon q^{-1/2}_{k_h-k_i}.\eqno\hbox{(2.13')}
$$
Replacing if necessary $\epsilon$ by $C_1\epsilon$ for a convenient constant $C_1$, only depending of $\alpha$ and $\beta$, we deduce from (2.12), (2.13), (2.14) and (2.13') that if we set $Q=x_{k_i,k_h}$, $R=y_{k_i,k_h}$, and $S=z_{k_i,k_h}$, then we have got: $$
\vert Q\alpha-R\vert\le\epsilon |Q|^{-1/2}, \ \ \ \vert Q\beta-S\vert\ll_{\alpha,\beta} |Q|^{-1/2},
\ \ \ \ Q\ne0,\ \ \ \ Q\equiv R\equiv0\ \ {\rm mod}\ D \eqno\hbox{(2.15)}$$
(where the implied constant in the inequality $\ll_{\alpha,\beta}$ only depends upon $\alpha$ and $\beta$).\par To make the proof complete, we must then find a triple of integers $(Q',R',S')$ satisfying condition (2.15), with $(Q',R',S')=1$. We shall use the following lemma:
\begin{Lemma} Let $(\alpha,\beta)$ be a pair of real numbers simultaneously badly approximable. 
Let $C_2$ be a positive real number. Suppose that a triple of integers $(Q,R,S)$, with $Q\ne0$, satisfies
 $$
\vert Q\alpha-R\vert\le|Q|^{-1/2}, \qquad \qquad\vert Q\beta-S\vert\le C_2|Q|^{-1/2}.$$ Denote by $g$ the {\rm g.c.d.}: $$g={\rm g.c.d.}(Q,R,S).$$ Then there exists a positive integer $G$, which depends only upon $(\alpha,\beta)$ and $C_2$, such that $g$ is a divisor of $G$.
 \end{Lemma} {\sl Proof}. Let $C$ be a positive constant such that for every triple of integers $(\hat q, \hat r,\hat s)$ with $\hat q>0$, we have $$\hat q^{1/2}\max\{\vert\hat q \alpha-\hat r\vert,\vert\hat q \beta-\hat s\vert\}\ge C,\eqno\hbox{(Bad2)}$$ and set
 $$ C_3=C^{-2/3}\max\{1,C_2^{2/3}\}.$$ Then we have:
$$0<g\le C_3.\eqno\hbox{(2.16)}$$
 Indeed, putting
$$Q=gq,\qquad R=gr, \qquad S=gs,$$  we have
$$ 
 |q|^{1/2}\vert q \alpha- r\vert\le g^{-3/2}, \ \ \qquad |q|^{1/2}\vert q\beta-s\vert\le C_2 g^{-3/2},
$$ and by (Bad2), we get (2.16). Then denote by $G$ the least common multiple of all the positive integers at most equal to $C_3$. It follows from (2.16) that $g$ divides $G$.\par
Now, to complete the proof of Theorem 1.3, we can suppose that $0<\epsilon\le1$, and we apply
Lemma 2.3 with the constant $C_2$ which is implied in (2.15). As this constant only depends on the pair $(\alpha,\beta)$, finally the integer $G$ depends only on $(\alpha,\beta)$. Replacing $D$ by $GD$, we can find a triple $(Q,R,S)$ of integers, with $Q\ne0$, such that
 $$
\vert Q\alpha-R\vert\le\epsilon |Q|^{-1/2}, \ \ \ \vert Q\beta-S\vert\le C_2|Q|^{-1/2},
\ \ \ \ Q\equiv R\equiv0\ {\rm mod}\ GD. \eqno\hbox{\rm(2.15A)}$$ Put $$g={\rm g.c.d.}(Q,R,S),$$
and 
$$Q=gQ',\qquad R=gR',\qquad S=gS'.$$ As $g$ divides $G$, we thus obtain integers $Q'\ne0$, $R'$, $S'$, such that
 $$
\vert Q'\alpha-R'\vert\le\epsilon|Q'|^{-1/2}, \ \ \ \vert Q'\beta-S'\vert\le C_2|Q'|^{-1/2},
\ \ \ \ Q'\equiv R'\equiv0\ {\rm mod}\ D,$$ with ${\rm g.c.d.}(Q',R',S')=1$. Thus Theorem 1.3 is proved.\par

 \section{Proof of Theorem 1.5} 
First we prove the following lemma: 
 \begin{Lemma} Let $p$ be a prime number, let $\Omega_p$ be the algebraic closure of ${\bf Q}_p$, and let $\zeta$ be an element of $\Omega_p$ of degree at most $3$ over ${\bf Q}_p$. Suppose that $$|\zeta|_p=1.$$ Then for each positive integer $\nu$,
$$|\zeta^{(p^6-1)p^\nu}-1|_p<p^{-\nu+1}.\eqno\hbox{\rm(3.1)}$$
\end{Lemma}
{\sl Proof}.
 First, note that if $u$ is an element of $\Omega_p$ with $|u|_p<1$, then
$$(1+u)^p=1+v$$ with $$|v|_p\le|u|_p\max\{1/p,|u|_p^{p-1}\}.\eqno\hbox{(3.2)}$$
Indeed, by the Newton formula, $$v=\sum_{k=1}^{p-1}{p\choose k} u^k+u^p$$ and for $0<k<p$, $$\left\vert{p\choose k}\right\vert_p=p^{-1},$$ so that $$\left\vert\sum_{k=1}^{p-1}{p\choose k} u^k\right\vert_p=p^{-1}|u|_p,$$ which proves (3.2). By induction, from (3.2), we deduce that, if $|u|_p\le p^{-1/(p-1)}$, then for each integer $\nu\ge0$
$$(1+u)^{p^\nu}=1+v_\nu$$ with $$|v_\nu|_p\le p^{-\nu}|u|_p.\eqno\hbox{(3.2')}$$
Now denote by $A_\zeta$ the ring of the elements $x\in{\bf Q}_p(\zeta)$ such that $|x|_p\le1$, and by $I_\zeta$ the maximal ideal of $A_\zeta$, i.e., $I_\zeta=\{x\in A_\zeta\ ;\ |x|_p<1\}$. As $\zeta$ has degree at most $3$ over ${\bf Q}_p$, the quotient field $A_\zeta/I_\zeta$ is ${\bf F}_p$, or ${\bf F}_{p^2}$, or ${\bf F}_{p^3}$. Also, for $x\in I_\zeta$, we have $|x|_p\le p^{-1/{\rm deg}\zeta}$. Since $|\zeta|_p=1$, we have in all the cases $$|\zeta^{p^6-1}-1|_p<1,$$ hence
 $$|\zeta^{p^6-1}-1|_p\le p^{-1/3}.$$ If $p\ge5$, we get by (3.2') $$|\zeta^{(p^6-1)p^\nu}-1|_p\le p^{-1/3-\nu},$$ and the result is proved in this case.
\par If $p=3$, and if $[{\bf Q}_p(\zeta):{\bf Q}_p]\le2$, then we have similarly for each $\nu\ge0$,
$$|\zeta^{(p^6-1)p^\nu}-1|_p\le p^{-1/2-\nu},$$ which proves Lemma 3.1 in this case. If $p=3$ with $[{\bf Q}_p(\zeta):{\bf Q}_p]=3$, we have just to look at the initial values; starting from
$$|\zeta^{p^6-1}-1|_p\le p^{-1/3},$$ which leads by (3.2) to $$|\zeta^{(p^6-1)p}-1|_p\le p^{-1},$$ we get by (3.2'):
$$|\zeta^{(p^6-1)p^\nu}-1|_p\le p^{-\nu}\qquad (\nu>0).$$
If $p=2$, and if $\zeta$ lies in ${\bf Q}_p$, then 
$$|\zeta^{p^6-1}-1|_p\le p^{-1},$$ and by (3.2'), we have for each $\nu\ge0$, $$|\zeta^{(p^6-1)p^\nu}-1|_p\le p^{-(\nu+1)}.$$ If $p=2$ and
$[{\bf Q}_p(\zeta):{\bf Q}_p]=2$, we have $$|\zeta^{p^6-1}-1|_p\le p^{-1/2},$$ then
$$|\zeta^{(p^6-1)p}-1|_p\le p^{-1},$$ and by (3.2'), for each $\nu>0$, $$|\zeta^{(p^6-1)p^\nu}-1|_p\le p^{-\nu}.$$
If $p=2$ with $[{\bf Q}_p(\zeta):{\bf Q}_p]=3$, then $$|\zeta^{p^6-1}-1|_p\le p^{-1/3},\qquad |\zeta^{(p^6-1)p}-1|_p\le p^{-2/3},\qquad |\zeta^{(p^6-1)p^2}-1|_p\le p^{-4/3},$$ thus by (3.2') $$
|\zeta^{(p^6-1)p^\nu}-1|_p\le p^{-\nu+2/3}\qquad (\nu\ge2).$$ Hence (3.1) is proved.\par Then we can prove the following lemma. Recall that in the introduction, we have defined the function $\psi$ by
$$\psi(D)=D\prod_{p|D}(p^6-1)$$ for each positive integer $D$.
\begin{Lemma}
Let $E$ be a cubic field and let $\theta$ be an algebraic integer in $E$ such that $${\rm Tr}(\theta)=0.$$ Let $\zeta$ be a unit of $E$. Then for each positive integer $D$, we have $${\rm Tr}(\theta\zeta^{\psi(D)})\equiv0\ \ {\rm mod}\ D.$$\end{Lemma}
{\sl Proof}.
Let $p$ be a prime number and let $\sigma_{j,p}$ ($j=0,1,2$) be the embeddings of $E$ into $\Omega_p$. Let $\nu$ be a positive integer. As ${\rm Tr}(\theta)=0$, we can write:
$${\rm Tr}(\theta\zeta^{(p^6-1)p^\nu})={\rm Tr}(\theta(\zeta^{(p^6-1)p^\nu}-1))=\sum_{0\le j\le2}\sigma_{j,p}(\theta)(\sigma_{j,p}(\zeta)^{(p^6-1)p^\nu}-1).$$ As $\theta$ is an algebraic integer, we have $|\sigma_{j,p}(\theta)|_p\le1$, and since $\zeta$ is a unit, we have $|\sigma_{j,p}(\zeta)|_p=1$.
Therefore, applying Lemma 3.1 to $\sigma_{j,p}(\theta)$, we get
$$|{\rm Tr}(\theta\zeta^{(p^6-1)p^\nu})|_p<p^{-\nu+1}.$$ Now ${\rm Tr}(\theta\zeta^{(p^6-1)p^\nu})$ lies in ${\bf Z}$, hence $|{\rm Tr}(\theta\zeta^{(p^6-1)p^\nu})|_p\in p^{-{\bf N}}$, and thus the above inequality leads to 
$$|{\rm Tr}(\theta\zeta^{(p^6-1)p^\nu})|_p\le p^{-\nu}.\eqno\hbox{(3.3)}$$ If $|D|_p=p^{-\nu}$ with $\nu>0$, then $\psi(D)=N(p^6-1)p^\nu$ where $N$ is a positive integer. As (3.3) also applies when replacing  $\zeta$ by $\zeta^N$, we have
$$|{\rm Tr}(\theta\zeta^{\psi(D)})|_p\le p^{-\nu}.\eqno\hbox{(3.3')}$$ Since this inequality holds for each prime divisor $p$ of $D$, we conclude that ${\rm Tr}(\theta\zeta^{\psi(D)})$ is divisible by $D$.\par
We can then prove Theorem 1.5. As the bilinear form $(x,y)\mapsto{\rm Tr}(xy)$ on $E$ is non degenerate, and since $(1,\alpha,\beta)$ is a basis of $E$ over ${\bf Q}$, there exists a number $\gamma_0\in E$ such that
$${\rm Tr}(\gamma_0)={\rm Tr}(\gamma_0\alpha)=0,\qquad{\rm Tr}(\gamma_0\beta)=1.$$ Thus, taking
$\gamma=\Gamma\gamma_0$ where $\Gamma$ is a convenient integer, we can find a number $\gamma\in E$
with $\gamma>0$, such that $\gamma$, $\gamma\alpha$ and $\gamma\beta$ are algebraic integers and 
$${\rm Tr}(\gamma)={\rm Tr}(\gamma\alpha)=0.\eqno\hbox{(3.4)}$$ In both the cases a) and b), the proof will be obtained with
$$Q={\rm Tr}(\gamma\zeta^{\psi(D)}),\qquad R={\rm Tr}(\gamma\alpha\zeta^{\psi(D)}),\qquad S={\rm Tr}(\gamma\beta\zeta^{\psi(D)})\eqno\hbox{(3.5)}$$ where $\zeta$ is a convenient unit of $E$.  Condition (3.5) implies in view of (3.4) and Lemma 3.2, that
$$Q\equiv R\equiv0\ {\rm mod}\ D.\eqno\hbox{(D)}$$
Let $\sigma_j$ ($j=1,2$) be the embeddings other than identity of $E$ into the complex field ${\bf C}$. In the case a), $\sigma_2=\overline{\sigma_1}$, and in the case b), $\sigma_1$ and $\sigma_2$ are both real. In order to choose the unit $\zeta$, we ask for some conditions. First, we require that there is a constant $C_4>1$, only depending upon $(\alpha,\beta)$, such that
$$\zeta\ge C_4,\qquad |\sigma_j(\gamma)||\sigma_j(\zeta)^{\psi(D)}|\le{\gamma\zeta^{\psi(D)}\over4}\qquad (j=1,2).\eqno\hbox{(3.6)}$$ Writing
$$Q=\gamma\zeta^D+\sigma_1(\gamma\zeta^D)+\sigma_2(\gamma\zeta^D),$$
we see that (3.6) implies that $Q$ will be a positive integer with
$${\gamma\zeta^{\psi(D)}\over2}\le Q\le{3\gamma\zeta^{\psi(D)}\over2}\cdot\eqno\hbox{(3.7)}$$  
Let us write
$$Q\alpha-R=\sigma_1(\gamma)(\alpha-\sigma_1(\alpha))\sigma_1(\zeta^{\psi(D)})+\sigma_2(\gamma)(\alpha-\sigma_2(\alpha))\sigma_2(\zeta^{\psi(D)})\eqno\hbox{(3.8)}$$
$$Q\beta-S=\sigma_1(\gamma)(\beta-\sigma_1(\beta))\sigma_1(\zeta^{\psi(D)})+\sigma_2(\gamma)(\beta-\sigma_2(\beta))\sigma_2(\zeta^{\psi(D)}),\eqno\hbox{(3.9)}$$ and note that it follows from (3.4) that 
$$Q\alpha-R=\sigma_1(\gamma)(\alpha-\sigma_1(\alpha))(\sigma_1(\zeta^{\psi(D)})-\sigma_2(\zeta^{\psi(D)})).\eqno\hbox{(3.8')}$$
We shall also require that
$$|\sigma_1(\zeta^{\psi(D)})-\sigma_2(\zeta^{\psi(D)})|\ll_{\alpha,\beta}\epsilon|\sigma_1(\zeta^{\psi(D)})|\eqno\hbox{(3.10)}$$ (the notation $\ll_{\alpha,\beta}$ means that the implied constant only depends upon $\alpha$ and $\beta$). In particular, taking $\epsilon$ small enough - in terms only depending on $(\alpha,\beta)$ -, we can then suppose that
$$|(\sigma_1(\zeta^{\psi(D)})-\sigma_2(\zeta^{\psi(D)})|\le{1\over2}|\sigma_1(\zeta^{\psi(D)})|,\eqno\hbox{(3.10')}$$ hence
$${1\over2}|\sigma_1(\zeta^{\psi(D)})|\le|\sigma_2(\zeta^{\psi(D)})|\le{3\over2}|\sigma_1(\zeta^{\psi(D)})|,$$ and since $|N_{E/{\bf Q}}(\zeta)|=|\zeta\sigma_1(\zeta)\sigma_2(\zeta)|=1$, we shall have
$$
|\sigma_j(\zeta^{\psi(D)})|\le\sqrt2|\zeta^{-\psi(D)/2}|\qquad(j=1,2).\eqno\hbox{(3.11)}$$ 
If we have got $\zeta\ge C_4>1$, then the second part of condition (3.6) will be ensured if
$$C_4^{3\psi(D)/2}\ge4\sqrt2{\max\{|\sigma_1(\gamma)|,\sigma_2(\gamma)|\}\over\gamma}\cdot\eqno\hbox{(3.6')}$$
Now note that in Theorem 1.5, we may replace $D$ by a multiple $D'$ of $D$ such that $D'/D$ is bounded (in terms only depending upon $\alpha$ and $\beta$). Accordingly, there is no loss of generality if we suppose $D\ge D_0$, where $D_0$ is an integer only depending upon the pair $(\alpha,\beta)$. Thus, as $\psi(D)\ge D$, we will obtain (3.6'), and (3.6) will follow.
From (3.7), (3.8'), (3.9), (3.10) and (3.11), we shall then deduce that $$|Q\alpha-R|\ll_{\alpha,\beta}\epsilon Q^{-1/2},\qquad |Q\beta-S|\ll_{\alpha,\beta} Q^{-1/2}.\eqno\hbox{(3.12)}$$
Suppose now that $E$ has a unique real embedding. Let $\eta$ be an irrational unit of $E$, replacing if necessary $\eta$ by $\eta^{\pm2}$, we can suppose that $\eta>1$. Let $\sigma$ and $\bar\sigma$ be not real embeddings of $E$ into ${\bf C}$. Since $N_{E/Q}(\eta)=\eta|\sigma(\eta)|^2=1$, we have $|\sigma(\eta)|=\eta^{-1/2}$
Let $\tau$ be a real number such that $$\sigma(\eta)=\eta^{-1/2}e^{i\tau}.$$ Thus for every positive integer $m$:
$$|\sigma(\eta^{m\psi(D)})-\bar\sigma(\eta^{m\psi(D)})|=\eta^{-m\psi(D)/2}|e^{2im\psi(D)\tau}-1|\le2\pi\eta^{-m\psi(D)/2}\Vert{m\psi(D)\tau\over\pi}\Vert$$  (where $\Vert.\Vert$ denotes the distance to the nearest integer). By Dirichlet's Theorem, if $1/\epsilon$ is an integer, we can choose an integer $m$ such that
$$0<m\le{1\over\epsilon},\qquad \Vert{m\psi(D)\tau\over\pi}\Vert<\epsilon.\eqno\hbox{(3.13)}$$ We thus get
 $$|\sigma(\eta^{m\psi(D)})-\bar\sigma(\eta^{m\psi(D)})|\le2\pi\epsilon\eta^{-m\psi(D)/2},\eqno\hbox{(3.10a)}$$ which ensures (3.10) by taking
 $\zeta=\eta^m\ge\eta>1$. We get also the first part of condition (3.6), hence in order to obtain the second part of (3.6), it is enough to suppose that $D$ is sufficiently large.
Thus the integers $Q$, $R$ and $S$ which are defined by (3.5) satisfy (3.12) together with the condition (D). Moreover
$\log\zeta^{\psi(D)}\ll_{\alpha,\beta}  m\psi(D)$, hence by (3.13)
$$\log Q\ll_{\alpha,\beta}  \epsilon^{-1}\psi(D),$$ which proves Theorem 1.5, a) without the condition $(Q,R,S)=1$. To add this last condition, Lemma 2.3 applies. Indeed it is well known that if $(1,\alpha,\beta)$ is a basis of a cubic field, then the pair $(\alpha,\beta)$ is simultaneously badly approximable. We replace $D$ by $GD$, where $G$ is defined as in Lemma 2.3, and we take a triple of integers $(Q,R,S)$ with $Q>0$, satisfying  (3.12), with $$Q\equiv R\equiv 0\ {\rm mod}\ GD,\qquad \log Q\ll_{\alpha,\beta}\epsilon^{-1}\psi(GD).$$ Define $$g={\rm g.c.d.}(Q,R,S)$$ and
$$Q'=\frac Qg,\qquad R'=\frac Rg,\qquad S'=\frac Sg\cdot$$ These integers satisfy (3.12), and since $g\ |\ G$, we get
$$Q'\equiv R'\equiv 0\ {\rm mod}\ D,\qquad {\rm g.c.d.}(Q',R',S')=1.$$ Moreover we have
$$\log Q'\le\log Q\ll_{\alpha,\beta}\epsilon^{-1} \psi(GD).\eqno\hbox{(1.13')}$$
Now $\psi(GD)\le\psi(G)\psi(D)$, and since $G$ only depends on $(\alpha,\beta)$, we get
$$\log Q'\ll_{\alpha,\beta}\epsilon^{-1} \psi(D),$$ which is (1.13). Thus Theorem 1.5, a) is proved.
\par
Suppose now that $E$ has two real embeddings $\sigma_j$ ($j=1,2$), other than identity. The group of units of $E$ has rank 2 \cite{Lang}, hence there exist multiplicatively independent units $\epsilon_1$ and $\epsilon_2$. Replacing $\epsilon_1$ and $\epsilon_2$ by $\epsilon_1^{\pm2}$ and $\epsilon_2^{\pm2}$, we can suppose that $\epsilon_1>1$, $\epsilon_2>1$ and $\sigma_j(\epsilon_\kappa)>0$ ($i=1,2,\ \kappa=1,2$). Given a positive integer $D$ and $0<\epsilon\le1$ such that $1/\epsilon$ is an integer, we construct a unit $\zeta=\epsilon_1^{m_1}\epsilon_2^{m_2}$ such that conditions (3.6) and (3.10) are satisfied. Note that $\sigma_1(\epsilon_j)\ne\sigma_2(\epsilon_j)$ for $j=1,2$, therefore
we can find by Dirichlet's Theorem integers $m_j\in{\bf Z}$ ($j=1,2$) such that
$$|m_1\left(\log\sigma_1(\epsilon_1)-\log\sigma_2(\epsilon_1)\right)+m_2\left(\log\sigma_1(\epsilon_2)-\log\sigma_2(\epsilon_2)\right)|
\le|\log\sigma_1(\epsilon_2)-\log\sigma_2(\epsilon_2)|{\epsilon\over\psi(D)}\eqno\hbox{(3.14)} $$ and $$0<|m_1|<{\psi(D)\over\epsilon}\cdot\eqno\hbox{(3.15)} $$
We get $$m_1\log\epsilon_1+m_2\log\epsilon_2=m_1\left(\log\epsilon_1-{\log\sigma_1(\epsilon_1)-\log\sigma_2(\epsilon_1)\over\log\sigma_1(\epsilon_2)-\log\sigma_2(\epsilon_2)}\log\epsilon_2\right)+{t\epsilon\over\psi(D)},$$ with $$|t|\le\log\epsilon_2.$$
As $\epsilon_j$ is a unit for $j=1,2$, we have $$\log\epsilon_j+\log\sigma_1(\epsilon_j)+\log\sigma_2(\epsilon_j)=0,$$ thus we get $$m_1\log\epsilon_1+m_2\log\epsilon_2=m_1(\log\epsilon_1-{2\log\sigma_1(\epsilon_1)+\log\epsilon_1\over2\log\sigma_1(\epsilon_2)+\log\epsilon_2}\log\epsilon_2)+{t\epsilon \over\psi(D)},$$ i.e., $$m_1\log\epsilon_1+m_2\log\epsilon_2=2{\log\epsilon_1\log\sigma_1(\epsilon_2)-\log\epsilon_2\log\sigma_1(\epsilon_1)\over2\log\sigma_1(\epsilon_2)+\log\epsilon_2}m_1+{t\epsilon\over\psi(D)} \cdot\eqno\hbox{(3.16)}$$
A classical result asserts that $$\log\epsilon_1\log\sigma_1(\epsilon_2)-\log\epsilon_2\log\sigma_1(\epsilon_1)\ne0$$ (see \cite{Lang}, p. 66). Then set $$M=\left\vert{\log\epsilon_1\log\sigma_1(\epsilon_2)-\log\epsilon_2\log\sigma_1(\epsilon_1)\over2\log\sigma_1(\epsilon_2)+\log\epsilon_2}\right\vert>0.$$ As before we can suppose that $D$ is large enough (in terms  only depending upon $\alpha$ and $\beta$) to get $$
{\log\epsilon_2\over\psi(D)}\le M.$$ Replacing if necessary the integers $m_j$ by their opposites, we have
$${\log\epsilon_1\log\sigma_1(\epsilon_2)-\log\epsilon_2\log\sigma_1(\epsilon_1)\over2\log\sigma_1(\epsilon_2)+\log\epsilon_2}m_1>0,$$
hence by (3.16), $$0<M\le M|m_1|\le m_1\log\epsilon_1+m_2\log\epsilon_2\le 3M|m_1|.\eqno\hbox{(3.17)}$$
Taking $\zeta=\epsilon_1^{m_1}\epsilon_2^{m_2}$, we consider $Q$, $R$ and $S$ defined by (3.5).
Condition (3.14) means that $$|\log\sigma_1(\zeta)-\log\sigma_2(\zeta)|\ll_{\alpha,\beta}{\epsilon\over\psi(D)}\cdot\eqno\hbox{(3.18)}$$
This inequality implies that 
$$\sigma_1(\zeta^{\psi(D)})\asymp_{\alpha,\beta}\sigma_2(\zeta^{\psi(D)}),$$ and thus we deduce from (3.18) that
$$|\sigma_1(\zeta^{\psi(D)})-\sigma_2(\zeta^{\psi(D)})|\ll_{\alpha,\beta} \epsilon\sigma_1(\zeta^{\psi(D)}),$$
which is the condition (3.10). By (3.17), if $D$ is large, condition (3.6) is satisfied with $C_4=e^{M}$. Hence condition (3.12) is satisfied together with condition (D). Moreover by (3.7), (3.15) and (3.17), $$\log Q\ll_{\alpha,\beta} |m_1|\psi(D)\ll_{\alpha,\beta} {\psi(D)^2\over\epsilon}\cdot$$ To achieve the proof of  Theorem 1.5, b), we must again ensure the condition ${\rm g.c.d.}(Q,R,S)=1$. We proceed as in the case a). We obtain then a triple $(Q',R',S')$ satisfying (3.12) with $$Q'\equiv R'\equiv0\ {\rm mod}\ D,\qquad {\rm g.c.d.}(Q',R',S')=1,$$ and
 $$\log Q'\ll_{\alpha,\beta} {\psi(GD)^2\over\epsilon},\eqno\hbox{(1.14')}$$ where $G$ is a positive integer only depending upon $(\alpha,\beta)$. As $\psi(GD)\ll_{\alpha,\beta}\psi(D)$, (1.14) immediately follows, thus Theorem 1.5, b) is proved. 
\section{Conclusion.} Unfortunately, the link between a pair $(\alpha,\beta)$ of real numbers and a sequence $(a_n,b_n,c_n)$ as defined in Lemma 1.2 is not clear. Given a finite set $\cal A$, if we consider the measure $\mu$ on ${\cal A}^{\bf N}$: $\mu=\bigotimes_{\bf N}\lambda$, where $\lambda$ is the measure on $\cal A$ such that for each set $X\subset\cal A$, $\lambda(X)={\rm card}(X)/\rm card\cal A$, then $\mu$-almost all infinite words formed with elements of $\cal A$ are recurrent. Thus, if we denote by $\cal B$ the set of simultaneously badly approximable pairs $(\alpha,\beta)$ for which there exists a word $(a_n,b_n,c_n)$ as defined in Lemma 1.2 which is recurrent, it is natural to think that ``many" pairs $(\alpha,\beta)$ satisfying (Bad2) are in $\cal B$. Nevertheless, we cannot prove this, even we cannot ensure that there exist in $\cal B$ other elements than pairs of numbers in the same cubic field. We cannot prove that $\cal B$ is uncountable. One can also ask whether $\cal B$ has full Hausdorff dimension. In the opposite sense, we do not know any argument which ensures that there exists simultaneouly badly  approximable pairs which are {\sl not} in $\cal B$. Notice that any pair of real numbers which does not satisfy (Bad2) satisfies the Littlewood conjecture, however we do not know whether such a pair satisfies (Ldiv). 
\par It would be interesting to extend to (Ldiv) some results which are known around the Littlewood conjecture. For instance, is it possible to prove that the set of exceptions for (Ldiv) has Hausdorff dimension zero? The same question may be restricted to simultaneously badly approximable pairs.
\par Our results can be extended to a field of formal power series over a finite ground field. It is well known that for an infinite ground field $K$, there exists a pair $(\alpha,\beta)\in K((T^{-1}))\times K((T^{-1}))$ which does not satisfy the Littlewood conjecture (\cite{DL}, \cite{Bk}). When the ground field $K$ is infinite, one can ask whether there exists a pair $(\alpha,\beta)\in K((T^{-1}))\times K((T^{-1})$ which satisfies the Littlewood conjecture and does not satisfy (Ldiv).
\par Notice that a counterexample to the {\sl mixed} Littlewood conjecture in ${\bf F}_3((T^{-1}))$ is given in \cite{Adal}.
\section{Acknowledgments} I thank the referee for his careful reading and his good questions which helped me to improve this paper. I thank also Yann Bugeaud for our interesting discussions.

  bernard.demathan@gmail.com
\end{document}